\documentclass[12pt,a4paper,draft]{amsart}
\usepackage{amsmath,amssymb,amsthm,bm,enumerate,longtable,color}
\usepackage[mathscr]{euscript}

\newcommand{\Aut}{\textnormal{Aut}}
\newcommand{\Sym}{\textnormal{Sym}}
\newcommand{\Rank}{\textnormal{Rank}}
\newcommand{\PG}{\textnormal{PG}}
\newcommand{\AGL}{\textnormal{AGL}}
\newcommand{\PairRank}{\textnormal{PairRank}}
\newcommand{\PP}{\mathscr{P}}
\newcommand{\D}{\mathscr{D}}
\newcommand{\F}{\mathbb{F}}
\newcommand{\B}{\mathscr{B}}

\newcommand{\C}{\mathscr{C}}
\newcommand{\OO}{\mathscr{O}}
\newcommand{\lan}{\langle}
\newcommand{\ra}{\rangle}

\def\la{\lambda}

\let\leq=\leqslant % redefine them
\let\geq=\geqslant

\newtheorem{theorem}{Theorem}[section]
\newtheorem{lemma}[theorem]{Lemma}
\newtheorem{corollary}[theorem]{Corollary}
\newtheorem{proposition}[theorem]{Proposition}

\theoremstyle{definition}
\newtheorem{definition}[theorem]{Definition}
\newtheorem{question}{Question}
\newtheorem{remark}[theorem]{Remark}
\newtheorem{example}[theorem]{Example}

\newtheorem{construction}[theorem]{Construction}

\title[Block-transitive point-imprimitive designs]{Delandtsheer--Doyen parameters for block-transitive point-imprimitive $2$-designs}
\author{Carmen Amarra, Alice Devillers, and Cheryl E. Praeger}
\thanks{The research in this paper forms part of Australian Research Council Discovery Project DP200100080.}
\date{\today}

\begin{document}

\begin{abstract}
Delandtsheer and Doyen bounded, in terms of the block size, the number of points of a point-imprimitive, block-transitive $2$-design. To do this they introduced two integer parameters $m,  n$, now called Delandtsheer--Doyen parameters, linking the block size with the parameters of an associated imprimitivity system on points.  We show that the Delandtsheer--Doyen parameters provide upper bounds on the permutation ranks of the groups induced on the imprimitivity system and on a class of the system.  We explore extreme cases where these bounds are attained, give a new construction for a family of designs achieving these bounds, and pose several open questions concerning the Delandtsheer--Doyen parameters.


\end{abstract}

\maketitle

\section{Introduction}\label{intro}

We study finite $2$-designs admitting a great deal of symmetry, and explore several extreme cases suggested by bounds on the so-called Delandtsheer-Doyen parameters. We consider \emph{$2$-$(v, k, \la)$ designs}: these are structures $\D = (\PP,\B)$ with two types of objects called \emph{points} (elements of $\PP$) and \emph{blocks} (elements of $\B$, sometimes called \emph{lines}). There are $v = |\PP|$ points, and we require that each block is a $k$-subset of $\PP$, and that each pair of distinct points lies in exactly $\la$ blocks. As a consequence of these conditions, the number $r$ of blocks containing a given point is also constant. The usual parameters associated with a $2$-$(v, k, \la)$ design are $v$, $k$, $r$, $\la$, and the number $b = |\B|$ of blocks, and we note that $v$, $k$, $\la$ determine $b$ and $r$ by
	\begin{equation}\label{eq1}
	\la(v-1) = r(k-1) \quad \mbox{and} \quad vr = bk, 	
	\end{equation}
see \cite[(5) on p. 57]{Dem}. \emph{Automorphisms} of $\D$ are permutations of $\PP$ that, in their induced action on $k$-subsets, leave invariant the block set $\B$. Thus the automorphism group $\Aut(\D)$ of the design is a subgroup of the symmetric group $\Sym(\PP)$. Automorphisms act on points, on blocks, and also on \emph{flags} (incident point-block pairs), and the following implications hold for transitivity of a subgroup $G \leq \Aut(\D)$ in these actions (the second implication follows from Block's Lemma, see \cite[(2.3.2)]{Dem}):
	\begin{center}
	flag-transitive \quad $\Rightarrow$ \quad block-transitive \quad $\Rightarrow$ \quad point-transitive.
	\end{center}

A celebrated result from 1961 of Higman and McLaughlin \cite{HM}, for $2$-designs with $\la = 1$, shows that flag-transitivity implies \emph{point-primitivity} for these designs, that is to say, the group does not leave invariant any \emph{non-trivial partition} of the point set (a partition with more than one class and class size at least $2$). Generalising this result, several other conditions on the parameters of a $2$-$(v,k,\la)$ design were given in the 1960s by Dembowski \cite[(2.3.7)]{Dem} and Kantor \cite[Theorems 4.7, 4.8]{Kan} under which flag-transitivity implies point-primitivity. Indeed point-primitivity is a desirable property for designs as it allows application of powerful theory for primitive permutation groups. Examples were known of block-transitive groups acting on $2$-designs which were not point-primitive, such as the ones we present in Example~\ref{ex:projpl}. Nevertheless it was hoped that (in some sense) most block-transitive groups on $2$-designs would be point-primitive, and this was shown to be the case by Delandtsheer and Doyen \cite{DD} in 1989. Their theorem, which we state below, implies that every block-transitive group will be point-primitive provided that the number of points is large enough, specifically $v > \left(\binom{k}{2} - 1\right)^2$ is sufficient. They consider a block-transitive group preserving a non-trivial point-partition, and they call an unordered pair of points an \emph{inner pair} if the two points lie in the same class of the partition, and an \emph{outer pair} if they lie in different classes. Since the group is block-transitive, the numbers of inner pairs and outer pairs in a block $B$ are constants, independent of the choice of $B\in\B$, and their sum is $\binom{k}{2}$.

\begin{theorem} \cite[Theorem]{DD} \label{theorem:DD}
Let $\D = (\PP, \B)$ be a $2$-$(v, k, \la)$ design, let $B \in \B$, and let $G \leq \Aut(\mathscr{D})$ be block-transitive. Suppose that $v = cd$ for some integers $c \geq 2$ and $d \geq 2$, and that $G$ leaves invariant a partition $\C$ of $\PP$ with $d$ classes, each of size $c$. Then there exist positive integers $m$ and $n$ such that
	\begin{equation} \label{eq:cd}
	c = \frac{\binom{k}{2} - n}{m} \quad \text{and} \quad
	d = \frac{\binom{k}{2} - m}{n}.
	\end{equation}
Moreover, $n$ is the number of inner pairs contained in $B$ and $mc$ is the number of outer pairs contained in $B$.
\end{theorem}

We call the integers $(m,n)$ the \emph{Delandtsheer--Doyen parameters} for $\D$ (relative to $G$ and $\C$). A major open question regarding these numbers is:

\begin{question} \label{q1}
Which Delandsheer--Doyen parameters $(m,n)$ are possible?
\end{question}

While these numbers have combinatorial significance, as given by Theorem~\ref{theorem:DD}, the purpose of this paper is to report on restrictions we discovered that the Delandtsheer--Doyen parameters place on the action of the group $G$. Let $K = G^\C$ denote the subgroup of $\Sym(\C)$ induced by $G$, and for $C \in \C$, let $H = G_C^C$ denote the subgroup of $\Sym(C)$ induced on $C$ by its setwise stabiliser $G_C$. By the Embedding Theorem \cite[Theorem 5.5]{PS}, we may assume that $G \leq H \wr K \leq \Sym(C) \wr \Sym(\C) \cong \Sym(c) \wr \Sym(d)$ in its imprimitive action on $\PP = \mathbb{Z}_c \times \mathbb{Z}_d$. For a transitive subgroup $X \leq \Sym(\Omega)$, the \emph{rank} of $X$ is the number $\Rank(X)$ of orbits in $\Omega$ of a point stabiliser $X_\alpha$ (for $\alpha\in\Omega$); and $\Rank(X)$ is also equal to the number of $X$-orbits in $\Omega\times\Omega$, see \cite[Lemma 2.28]{PS}. Similarly  we denote by $\PairRank(X)$ the number of $X$-orbits on the unordered pairs of distinct points from $\Omega$,
 and it is not difficult to see that $(\Rank(X)-1)/2 \leq \PairRank(X)\leq \Rank(X)-1$.  A summary of the major restrictions we obtain on the Delandsheer--Doyen parameters is given by the following theorem.

\begin{theorem} \label{theorem:rankbounds}
Let $\D$, $G$, $\C$, $c$, $d$, $m$, and $n$ be as in Theorem $\ref{theorem:DD}$. Let $C \in \C$ and $H = G_C^C$, $K = G^\C$. Then
	\[ \frac{\Rank(H)-1}{2} \leq \PairRank(H)\leq n,
	\quad \mbox{and}\quad
	\frac{\Rank(K)-1}{2} \leq \PairRank(K)\leq m. \]
\end{theorem}

Theorem~\ref{theorem:rankbounds} follows immediately from Proposition~\ref{prop:HK}, which contains additional detailed information about the permutation actions of the groups $H$ and $K$. Refining Question~\ref{q1} we might ask:

\begin{question} \label{q2}
For what values of $(m,n)$ can we have $\Rank(H) = 2n + 1$ and $\Rank(K) = 2m + 1$? For what values of $(m,n)$ can we have $\PairRank(H) = n$ and $\PairRank(K) = m$?
\end{question}

Some rather incomplete answers to Question~\ref{q2} may be deduced from certain results and examples of block-transitive designs in \cite{CP, OPP, NNOPP}. Examining these from the point of view of  the Delandsheer--Doyen parameters, we obtain the next result.

\begin{theorem} \label{theorem:exist}
	\begin{enumerate}[(a)]
	
	\item For any $N > 0$, there exists an example in Theorem $\ref{theorem:rankbounds}$ with $\Rank(H) = 2n + 1$, $\Rank(K) = 2m + 1$, and with both $n, m > N$.
	
	\item For any $k \geq 3$, there exist examples in Theorem $\ref{theorem:rankbounds}$ with $m=n=\PairRank(H) = \PairRank(K) = 1$, but with $\Rank(H) = \Rank(K) = 2$.
	
	\item There exists a $2$-$(cd, k, \la)$ design $\D$ (for some $\la$) in Theorem $\ref{theorem:rankbounds}$ with $n = m = 1$, and $\max\{ \Rank(H), \Rank(K) \} = 3$, if and only if $k = 8$ and $c = d = 27$. Moreover, in the case  $\la = 1$, there are up to isomorphism exactly $467$ such designs.
	
	\end{enumerate}
\end{theorem}

We note that additional examples may be constructed of block-transitive, point-imprimit\-ive $2$-$(729, 8, \la)$ designs in Theorem~\ref{theorem:exist}(c) using the technique from \cite[Proposition 1.1]{CP} applied to a subgroup of $\Sym(c) \wr \Sym(c)$ properly containing the group $G$ of Theorem~\ref{theorem:rankbounds}. Such examples will have larger values of $\la$. Theorem~\ref{theorem:exist} will be proved in Section~\ref{sec:existence}.

In the final Section~\ref{sec:construction} we present a new design construction that yields additional pairs $(m,n)$ with the second property requested in Question~\ref{q2}. The construction is different from, but was inspired by, the design construction in \cite[Proposition 2.2]{CP}. The construction in Section~\ref{sec:construction} requires integer pairs $[n,c]$ with the following property.

\begin{definition} \label{def:nc}
An integer pair $[n,c]$ is said to be \emph{useful} if the following two conditions hold:
	\begin{enumerate}
	\item $n \geq 2$ and $c$ is a prime power such that $c \equiv 1 \pmod{2n}$; and
	\item $c + n = \binom{k}{2}$ for some integer $k \geq 2n$.
	\end{enumerate}
\end{definition}

We show in Lemma~\ref{lemma:kbounds} that for useful pairs $[n,c]$, the value of $k$ satisfies $2n + 2 \leq k \leq n + d$, where $d = 1 + (c-1)/n$. Table~\ref{table:nc} gives a list of all useful pairs $[n,c]$ such that $n \leq 20$ and $c \leq 1300$, together with the corresponding values of $k$ and $d$. In addition, for $n\in\{11, 13, 16, 18\}$, the table contains the parameters $c, k, d$, for the smallest value of $c$ such that $[n, c]$ is useful. We note that the only integers $n$ in the range $2\leq n\leq 20$ which do not appear in the table are $n\in\{6, 10, 15\}$, and we prove in Lemma~\ref{lemma:nwithnousefulpair} that for these three values of $n$ there is no  $c$ such that $[n,c]$ is a useful pair.

%Old version of this table after \end{document}
\begin{table}[h]
\label{table:nc}
\begin{tabular}{crrr||crrr||crrr||crrr}
$n$ & $c$ &$k$ & $d$ &
	$n$ & $c$ &$k$ & $d$ &
	$n$ & $c$ &$k$ & $d$ &
	$n$ & $c$ &$k$ & $d$ \\
\hline \hline
$2$ &  $13$ & $6$ & $7$ &
	$3$ &  $25$ & $8$ & $9$ &
	$5$& $1171$& $49$& $235$&
	$12$ &  $313$ & $26$ & $27$ \\
$2$ &  $53$ & $11$ & $27$ &
	$4$ &  $41$ & $10$ & $11$ &
	$7$ &  $113$ & $16$ & $17$&
	$13$&$2003$&$64$&$155$\\
$2$ &  $89$ & $14$ & $45$ &
	$4$ &  $857$ & $42$ & $215$ &
	$7$ &  $659$ & $37$ & $95$ &
	 $14$ &  $421$ & $30$ & $31$ \\
$2$ &  $169$ & $19$ & $85$ &
	$5$ &  $61$ & $12$ & $13$ &
	$8$& $1217$& $50$& $153$ &
	 $16$&$25409$&$226$&$1589$\\
$2$ &  $229$ & $22$ & $115$ &
	$5$ &  $131$ & $17$ & $27$ &
	$9$ &  $181$ & $20$ & $21$ &
	$17$ &  $613$ & $36$ & $37$ \\ 
$2$ &  $349$ & $27$ & $175$ &
	$5$ &  $271$ & $24$ & $55$ &
	$9$ &  $397$ & $29$ & $45$ &
	$18$ &$1693$   &$59$  &$95$ \\ 
$2$ &  $433$ & $30$ & $217$ &
	$5$ &  $401$ & $29$ & $81$ &
	$9$ &  $487$ & $32$ & $55$ &
	$19$ &  $761$ & $40$ & $41$ \\
$2$ &  $593$ & $35$ & $297$ &
	$5$ &  $491$ & $32$ & $99$ &
	$9$ &  $811$ & $41$ & $91$ &
	$20$ &  $841$ & $42$ & $43$ \\  
	
$2$ &  $701$ & $38$ & $351$ &
	$5$ &  $661$ & $37$ & $133$&
	$9$ &  $937$ & $44$ & $105$ &
	&   &  & \\
$2$& $1033$& $46$& $517$ & 
$5$ &  $941$ & $44$ & $189$  &
 $11$&$2069$&$65$&$189$&
 &   &  &  \\
 \hline
\end{tabular}
\caption{Examples of useful pairs $[n,c]$, together with the  values for $k$ and $d$.}
\end{table}

We are principally interested in which values of $n$ are possible since, in our design Construction~\ref{con:newdesign} based on a useful pair $[n,c]$, the Delandtsheer--Doyen parameters turn out to be $(1,n)$.  Moreover, these designs also  satisfy the bounds $\PairRank(H) = n$ and $\PairRank(K) = 1$ in Question~\ref{q2} (see Theorem~\ref{theorem:newdesign}). The following theorem is an immediate consequence of Theorem~\ref{theorem:newdesign}. We note that an expression for the parameter $\la$ is determined in Remark~\ref{rem:newdesign-la}.

%Our interest in these pairs stems from the fact that, for each useful pair $[n,c]$, Construction~\ref{con:newdesign} produces a block-transitive, point-imprimitive $2$-design with Delandtsheer--Doyen parameters $(1,n)$ satisfying the bounds $\PairRank(H) = n$ and $\PairRank(K) = 1$ in Question~\ref{q2} (see Theorem~\ref{theorem:newdesign}). The following theorem is an immediate consequence. We note that an expression for the parameter $\la$ is determined in Remark~\ref{rem:newdesign-la}.

\begin{theorem} \label{theorem:existence2}
Suppose that $[n,c]$ is useful with $c + n = \binom{k}{2}$, and let $d = 1 + (c-1)/n$. Then there exists a $2$-$(cd, k, \la)$ design (for some $\la$) admitting a block-transitive, point-imprimitive group $H \wr K$ with $H \leq \Sym(c)$ and $K = \Sym(d)$, and with Delandtsheer--Doyen parameters $(1,n)$ such that $\Rank(H) = \PairRank(H) + 1 = n + 1$ and $\Rank(K) = \PairRank(K) + 1 = 1 + 1 = 2$.
\end{theorem}

Apart from the examples given to prove Theorems~\ref{theorem:exist}, \ref{theorem:existence2}, and~\ref{theorem:newdesign}, Questions~\ref{q1} and~\ref{q2} are in general wide open, and we would be very interested in knowing more general answers. 
In particular, we note that Theorem~\ref{theorem:existence2} does not produce designs  
with Delandtsheer--Doyen parameters $(1,n)$ for  $n=6, 10$ or $15$, since by 
Lemma~\ref{lemma:nwithnousefulpair}, there are no useful pairs $[n,c]$ for these values of $n$. Nevertheless there might be alternative constructions with Delandtsheer--Doyen parameters $(1,n)$ for such $n$.

\begin{question} \label{q3}
%Is it true that, for each integer $n \geq 2$, there exists at least one prime power $c$ such that $[n,c]$ is useful? 
For which values of $n$ do examples exist with Delandtsheer--Doyen parameters $(1,n)$, and with $(\PairRank(K), \PairRank(H))=(1,n)$? Are there examples for all $n$?
\end{question}

In particular, it would be good to know for which values of $n$ there exists at least one useful pair $[n,c]$.
We finish this introductory section with some commentary on number theoretic questions related to the existence of useful pairs.

%We prove in Lemma \ref{lemma:nwithnousefulpair} that no useful pair exists for $n=6,10,15$ however there is exactly one tuple  $[n,c,k,d]$ for each satisfying part (1) and the first part of (2) in Definition~\ref{def:nc}, they  are $[6, 49, 11, 9]$ and $[10, 81, 14, 9]$ and $[15, 121, 17, 9]$, but unfortunately $k < 2n$ in these three cases.
% We notice in particular that all values of $n$ up to $12$ arise except for $6$ and $10$. In fact further searching revealed that, for $n \in \{6, 10\}$, there is no value of $c \leq 10^6$ for which $[n,c]$ is useful, and the only tuples $[n,c,k,d]$ in this range satisfying part (1) and the first part of (2) in Definition~\ref{def:nc} are $[6, 49, 11, 9]$ and $[10, 81, 14, 9]$, but unfortunately $k < 2n$ in these two cases. Moreover, for $n = 6$, since $c \equiv 1 \pmod{12}$, we must have $k \equiv 11$ or $14 \pmod{24}$ and a search for such $k \leq 24 \times 10^7 + 14$ yielded no possibilities with $c = \binom{k}{2} - 6$ a prime power, apart for the one mentioned with $c = 49$. \emph{We wonder if perhaps there is no $c$ such that $[6,c]$ is useful?}
%
%

\subsection{Useful pairs and a conjecture from number theory.}\label{sub:useful}
By Dirichlet's Theorem on arithmetic progressions (see \cite[Chapter VIII.1]{SMC}), for any positive integer $n$, there exist infinitely many primes $c$ such that $c \equiv 1 \pmod{2n}$. However it is unclear how many such pairs $[n,c]$ sum to a triangular number $\binom{k}{2}$. In the light of Dirichlet's Theorem and the relatively large number of useful pairs of the form $[2,c]$ we found 
with $c < 1300$, we ask:

\begin{question}\label{q3a}
Are there  infinitely many prime powers $c$ such that $[2,c]$ is a useful pair?  
\end{question}

The conditions for $[n,c]$ to be useful imply in particular that, if $k\equiv r\pmod{4n}$, then $\binom{k}{2}\equiv \binom{r}{2}\equiv n+1\pmod{2n}$. Thus, for fixed positive integers $n, r$ such that $r<4n$ and $\binom{r}{2}\equiv n+1\pmod{2n}$, we are concerned with integers $k$ of the form $k=4nb+r$ for some positive integer $b$. The quantity $c=\binom{k}{2}-n$ is then a quadratic polynomial $f(b)$ with integer coefficients determined by $r$ and $n$, and we would like $f(b)$ to be a prime power. For example, if $n=2$ and $r=3$, we get $f(b)=32b^2+20b+1$.  
%
%If we set $k=4nb+r$ for a fixed $r$ such that $ \binom{k}{2}-n\equiv 1 \pmod {2n}$, we get  $ \binom{k}{2}-n$ as a quadratic polynomial $f(b)$ with integer coefficients. 
 
Let $f(b)$ be a polynomial  in the variable $b$ with integer coefficients. We claim that, if the sequence $f(1),f(2),f(3),\ldots $\quad contains infinitely many primes, then the following three conditions must hold: 
 
 \begin{tabular}{cl}
    (i)  &  the leading coefficient is positive;  \\
    (ii) & the polynomial is irreducible over the integers; and \\
     (iii) &   $\gcd(f(1),f(2),f(3),\ldots)=1$.
 \end{tabular}
 
 \noindent
 Condition (i) must hold as otherwise $f(n)$ takes only finitely many positive values. Condition (ii) must hold since $f(n)=g(n)h(n)$ being prime implies $g(n)=\pm 1$ or $h(n)=\pm 1$, but this can only happen for finitely many values of $n$. Condition (iii) must hold since otherwise there is a prime $p$ dividing $f(n)$ for all $n\geq1$, and then the only way $f(n)$ can be prime is if $f(n)=p$, which can only happen for a finite number of integers $n$.
 
 These three conditions are satisfied, for example, by the polynomial  $f(b)=32b^2+20b+1$ we mentioned above. The Russian mathematician Viktor Bunyakovsky (or Bouniakowsky) conjectured in 1857 that these  three conditions are also sufficient,  see \cite{B}. %
 To the best of our knowledge his conjecture is still open, apart from the degree $1$ case which is Dirichlet's Theorem. If the Bunyakovsky Conjecture were true then, for example,  $f(b)=32b^2+20b+1$ would be prime for infinitely many integers $b$, and hence there would be infinitely many useful pairs $[2,c]$ with $c$ a prime, answering Question~\ref{q3a}. Thus in general we ask:

\begin{question} \label{q3b}
For which integers $n$ do there exist infinitely many useful pairs $[n,c]$? 
\end{question}

%%%%%%%%%%%%%%%%%%%%%%%%%%%%%%%%%%%%%%%%%%%%%%%%%%%%%%%%%%%%%%%%%%%

\section{Permutation group concepts} \label{sec:perm}

Let $X$ be a transitive permutation group on a set $\Omega$. An \emph{$X$-orbital} is an $X$-orbit in $\Omega \times \Omega$. Clearly, $\{ (\alpha,\alpha) \ | \ \alpha \in \Omega \}$ is an orbital and is called the \emph{trivial orbital}; all other orbitals are said to be non-trivial.

For any $X$-orbital $\Delta$ and any $\alpha \in \Omega$, the set $\Delta(\alpha) = \{ \beta \ | \ (\alpha,\beta) \in \Delta \}$ is an $X_\alpha$-orbit, and is called a \emph{suborbit} of $X$. The set of $X$-orbitals is in one-to-one correspondence with the set of all $X_\alpha$-orbits in $\Omega$, such that the orbital $\Delta$ corresponds to the $X_\alpha$-orbit $\Delta(\alpha)$. In particular, the trivial orbital corresponds to the trivial suborbit $\{\alpha\}$.

The cardinality $|\Delta(\alpha)|$ is a \emph{subdegree} of $X$, and the number of $X$-orbitals (including the trivial orbital) is the \emph{rank} of $X$, denoted $\Rank(X)$.

For each $X$-orbital $\Delta$, the set $\Delta^* = \{ (\beta,\alpha) \ | \ (\alpha,\beta) \in \Delta \}$ is also an $X$-orbital, called the \emph{paired orbital} of $\Delta$. If $\Delta = \Delta^*$, then $\Delta$ is said to be \emph{self-paired}. For any $\alpha \in \Omega$, the set $\Delta(\alpha) \cup \Delta^*(\alpha)$ is therefore either a single suborbit, or the union of two suborbits of equal lengths. We call the cardinality
	\[ u_\Delta := \big|\Delta(\alpha) \cup \Delta^*(\alpha)\big| \]
the \emph{symmetrised subdegree} corresponding to $\Delta$ (or to $\Delta^*$). Note that $u_\Delta = \delta_\Delta |\Delta(\alpha)|$ where $\delta_\Delta = 1$ or $2$ according as $\Delta(\alpha)$ is self-paired or not. Let $\OO_X$ denote the set of all $\{\Delta, \Delta^*\}$, where $\Delta$ is a non-trivial $X$-orbital. Then $\OO_X$ is in one-to-one correspondence with the set of $X$-orbits on the unordered pairs of distinct points from $\Omega$, and we call $|\OO_X|$ the \emph{pair-rank} of $X$, denoted $\PairRank(X)$. It follows from the definition that 
	\begin{equation} \label{eq:pr}
	\PairRank(X) + 1 \leq \Rank(X) \leq 2\,\PairRank(X) + 1.
	\end{equation}

%%%%%%%%%%%%%%%%%%%%%%%%%%%%%%%%%%%%%%%%%%%%%%%%%%%%%%%%%%%%%%%%%%%

\section{Proof of Theorem \ref{theorem:rankbounds}}

Let $\D = (\PP,\B)$ be a $2$-$(v, k, \la)$ design, with $v = cd$ for some integers $c \geq 2$ and $d \geq 2$. Suppose that $G \leq \Aut(\D)$ is transitive on the block set $\B$ and leaves invariant a non-trivial partition $\C$ of $\PP$ with $d$ classes $C_1, \ldots, C_d$, each of size $c$.
The following lemma establishes useful identities between the parameters.
\begin{lemma}\label{counting}
Let $\D$, $\C$, $c,d$ be as above. Let $m$ and $n$ satisfy Equation \eqref{eq:cd}. Then the following identities hold.
	
\begin{enumerate}[(a)]
\item $cd-1= \binom{k}{2} \cdot \frac{c-1}{n}	= \binom{k}{2} \cdot \frac{d-1}{m}$;
\item $m(c-1)=n(d-1)$;
    \item The number of blocks is $|\B| =\frac{ cd(c-1)\la}{2n}=\frac{ cd(d-1)\la}{2m}$.
\end{enumerate}
\end{lemma}
\begin{proof}
\begin{enumerate}[(a)]
\item 
By Equation \eqref{eq:cd}, we have
	\[ cd - 1
	= \frac{\binom{k}{2} - n}{m} \cdot \frac{\binom{k}{2} - m}{n} - 1
	= \frac{\binom{k}{2}\left(\binom{k}{2} - n - m\right)}{nm}
	= \binom{k}{2} \cdot \frac{c-1}{n}
	= \binom{k}{2} \cdot \frac{d-1}{m}. \]
	
\item Part (b) follows immediately from Part (a).
\item  Note that $|\B| = v(v-1)\la/(k(k-1)) = cd(cd-1)\la/(k(k-1))$ 
by \eqref{eq1}, and
hence $|\B| = cd(c-1)\la/(2n)=cd(d-1)\la/(2m)$ using Part (a).
\end{enumerate}
\end{proof}

Let $K = G^\C$ denote the induced action of $G$ on the set $\C = \{C_1, \ldots, C_d\}$ of imprimitivity classes, and for $C \in \C$ let $H = G_C^C$ denote the induced action on $C$ of the setwise stabiliser $G_C$. Then by the Embedding Theorem for transitive permutation groups, we may assume that $G \leq H \wr K \leq \Sym(C) \wr \Sym(\C) \cong \Sym(c) \wr \Sym(d)$, \cite[Theorem 5.5]{PS}.

Let $X = H^d = H_1 \times \cdots \times H_d$ be the base group of the wreath product $H \wr K$, such that for each $i \in \{1, \ldots, d\}$, $H_i \cong H$ and $H_i$ induces $H$ on $C_i$ and fixes all other classes pointwise. Let $\Sigma$ be a non-trivial $H$-orbital. Then for each $i$, there is a corresponding $H_i$-orbital $\Sigma_i$ for the action of $H_i$ on $C_i$.

\begin{proposition} \label{prop:HK}
Let $\D$ and $G$ be as above and let $B \in \B$. Let $m$ and $n$ be as in Theorem $\ref{theorem:DD}$.
	\begin{enumerate}[(a)]
	\item For a non-trivial orbital $\Delta$ of $K$ with symmetrised subdegree $u_\Delta$, the number $nu_\Delta/(c-1)$ is an integer and there are exactly $cn u_\Delta/(c-1)$ pairs $\{\alpha, \beta\}$ in $B$ such that $\alpha \in C_i$ and $\beta \in C_j$ for some $(C_i,C_j) \in \Delta$. 	Moreover, $(\Rank(K)-1)/2 \leq \PairRank(K) \leq m$.
	\item For a non-trivial orbital $\Sigma$ of $H$ with symmetrised subdegree $u_\Sigma$, the number $n u_\Sigma/(c-1)$ is an integer and is equal to the number of pairs $\{\alpha, \beta\}$ in $B$ such that $(\alpha,\beta) \in \Sigma_i$ for some $i \in \{1, \ldots, d\}$. Moreover, $(\Rank(H)-1)/2 \leq \PairRank(H) \leq n$.
	\end{enumerate}
\end{proposition}

Part of (a) is proved in \cite[Lemma 2.1]{OPP}, but with different notation so we give brief details here (note that our parameter $u_\Delta$ is equal to the expression $2u/\delta$ in that reference).

\begin{proof}
We first prove part (a). For $\Delta$ a non-trivial orbital of $K$, let
	\[ 
	\mathscr{S}(\Delta) = \big\{ \{\alpha,\beta\} \ | \ \alpha \in C_i, \ \beta \in C_j \ \text{for some} \ (C_i,C_j) \in \Delta \big\}. 
	\]
The number of choices of $(C_i,C_j)$ is $|\Delta| = d|\Delta(C_i)| = du_\Delta/\delta_\Delta$, and for each choice, there are $c^2$ pairs $\{\alpha,\beta\} \in \mathscr{S}(\Delta)$ with $\alpha \in C_i, \beta \in C_j$. If $\Delta = \Delta^*$, that is, if $\delta_\Delta = 1$, then we have counted each unordered pair $\{\alpha,\beta\} \in \mathscr{S}(\Delta)$ twice so $|\mathscr{S}(\Delta)| = c^2du_\Delta/2$, while if $\Delta \ne \Delta^*$, that is, if $\delta_\Delta = 2$, then there is no double counting, and $|\mathscr{S}(\Delta)| = c^2du_\Delta/\delta_\Delta = c^2du_\Delta/2$. Hence $|\mathscr{S}(\Delta)| = c^2du_\Delta/2$ in either case.

Since $G$ leaves $\mathscr{S}(\Delta)$ invariant and is transitive on $\B$, each block $B$ contains the same number of pairs from $\mathscr{S}(\Delta)$, say $n_\Delta$ pairs. Thus, counting pairs $(\{\alpha,\beta\}, B)$ with $\{\alpha,\beta\}  \in \mathscr{S}(\Delta)$, $B \in \B$, and $\{\alpha,\beta\} \subseteq B$, we obtain $|\B| n_\Delta = \la |\mathscr{S}(\Delta)| = \la c^2du_\Delta/2$.  Since $|\B| = cd(c-1)\la/(2n)$ by Lemma \ref{counting}(c), it follows that $n_\Delta = cnu_\Delta/(c-1)$. In particular $nu_\Delta/(c-1)$ is an integer.

Note that $\mathscr{S}(\Delta^*) = \mathscr{S}(\Delta)$, and that each outer pair in $B$ lies in exactly one of the sets $\mathscr{S}(\Delta)$. Thus the number of outer pairs in $B$, namely $mc$ by Theorem~\ref{theorem:DD}, is equal to the sum of the $n_\Delta$ over the set $\OO_K$ of all pairs $\{\Delta,\Delta^*\}$ of non-trivial $K$-orbitals; that is,
	\[ 
	mc = \sum_{\{\Delta,\Delta^*\} \in \OO_K} \frac{cnu_\Delta}{c-1}. 
	\]
Since each $nu_\Delta/(c-1)$ is an integer, each term of the summation above is a positive integer multiple of $c$, Hence $m$ is greater than or equal to $|\OO_K|$ which, as we noted in Section~\ref{sec:perm}, is equal to $\PairRank(K)$. Thus, using \eqref{eq:pr},
	\[ 
	m \geq |\OO_K| = \PairRank(K) \geq (\Rank(K) - 1)/2. 
	\]

Now we prove part (b). Let $\Sigma$ be a non-trivial orbital of $H$. Note that $|\Sigma| = c|\Sigma(\alpha)| = cu_\Sigma/\delta_\Sigma$ and that, for each $i\leq d$, there is a corresponding $H_i$-orbital $\Sigma_i$. Let
	\[ 
	\mathscr{S}'(\Sigma) = \big\{ \{\alpha,\beta\} \ | \ (\alpha,\beta) \in \Sigma_i \ \text{for some} \ i \in \{1, \ldots, d\} \big\}. 
	\]
Since $G$ is transitive on $\mathscr{C}$, the set $\mathscr{S}'(\Sigma)$ contains equally many pairs from each class in $\mathscr{C}$. So for a fixed class $C \in \mathscr{C}$, and viewing $\Sigma$ as an $H$-orbital in $C\times C$,
	\begin{align*}
	|\mathscr{S}'(\Sigma)|
	&= d \cdot \big|\big\{ \{\alpha,\beta\} \ | \ \alpha, \beta \in C; \ \{\alpha,\beta\} \in \mathscr{S}'(\Sigma) \big\}\big| \\
	&= d \cdot \big|\big\{ \{\alpha,\beta\} \ | \ (\alpha,\beta) \in \Sigma \big\}\big|.
	\end{align*}
If $\Sigma = \Sigma^*$ then $\delta_\Sigma = 1$ and each unordered pair $\{\alpha,\beta\} \in \mathscr{S}'(\Sigma)$ from $C$ is counted twice (since both $(\alpha,\beta)$ and $(\beta,\alpha)$ lie in $\Sigma$), so $|\mathscr{S}'(\Sigma)| = d|\Sigma|/2 = dcu_\Sigma/2$. On the other hand, if $\Sigma \ne \Sigma^*$, then $\delta_\Sigma = 2$ and $|\mathscr{S}'(\Sigma)| = d|\Sigma| = dcu_\Sigma/2$. Hence in both cases $|\mathscr{S}'(\Sigma)| = dcu_\Sigma/2$. Since $G$ leaves $\mathscr{S}'(\Sigma)$ invariant and is transitive on $\mathscr{B}$, each block $B$ contains the same number of pairs from $\mathscr{S}'(\Sigma)$, say $n_\Sigma$ pairs. Thus, counting pairs $(\{\alpha,\beta\},B)$ with $\{\alpha,\beta\} \in \mathscr{S}'(\Sigma)$, $B \in \B$, and $\{\alpha,\beta\} \subseteq B$, we obtain $|\mathscr{B}|n_\Sigma = \la |\mathscr{S}'(\Sigma)| = \la dcu_\Sigma/2$. 
%Recall from above that $|\B| = cd(c-1)\la/(2n)$, so 
%Next we count in two ways the number of pairs $(p,B)$ such that $p$ is an inner pair of points, $B \in \B$, and $p \subseteq B$. This gives us $|\B|n = d\binom{c}{2}\la$, and hence $|\B| = d\binom{c}{2}\la/n$. Thus, using the expression above for $|\B|n_\Sigma$, we obtain
Using Lemma \ref{counting}(c) we obtain that $n_\Sigma=nu_\Sigma/(c-1).$
%	\[ n_\Sigma
%	= \frac{\la dcu_\Sigma}{2} \cdot \frac{1}{|\mathscr{B}|}
%	= \frac{\la dcu_\Sigma}{2} \cdot \frac{2n}{cd(c-1)\la}
%	= \frac{n u_\Sigma}{c-1}. \]
In particular, $nu_\Sigma/(c-1)$ is a positive integer. Note that $\mathscr{S}'(\Sigma) = \mathscr{S}'(\Sigma^*)$ for any non-trivial $H$-orbital $\Sigma$, and that each inner pair in $B$ lies in exactly one of the sets $\mathscr{S}'(\Sigma)$. Thus the number of inner pairs in $B$ is equal to the sum of the numbers $n_\Sigma$ over the set $\OO_H$ of all pairs $\{\Sigma,\Sigma^*\}$ of non-trivial $H$-orbitals. Hence
	\begin{equation} \label{eq:n}
	n = \sum_{\{\Sigma,\Sigma^*\} \in \OO_H} \frac{n u_\Sigma}{c-1}.
	\end{equation}
So, using \eqref{eq:pr}, $n \geq |\OO_H| = \PairRank(H) \geq (\Rank(H) - 1)/2$, which completes the proof.
\end{proof}

Proposition~\ref{prop:HK} has the following corollary. It is easy to prove: the condition $m = 1$ implies by Proposition~\ref{prop:HK}(a) that $\PairRank(K) = 1$, that is to say, $K$ is transitive on unordered pairs of distinct classes of $\C$. This means in particular that $K$ is primitive on $\C$, see \cite[Lemma 2.30]{PS}. Similarly, $n = 1$ implies that $H$ is primitive on $C$. Part (b) of this corollary was proved also in \cite[Lemma 2.3]{OPP}.

\begin{corollary}
Let $\D$, $G$, $H$, $K$, $\C$, $C$, $m$, and $n$ be as in Proposition~$\ref{prop:HK}$.
	\begin{enumerate}[(a)]
	\item If $m = 1$ then $K$ is primitive on $\C$.
	\item If $n = 1$ then $H$ is primitive on $C$.
	\end{enumerate}
\end{corollary}

Our last result of this section looks at cases where the upper bound on $\Rank(K)$ or $\Rank(H)$ is sharp. A transitive permutation group is {\it $3/2$-transitive} if all its non-trivial suborbits have the same size. 

\begin{lemma} \label{lemma:maxrank}
	\begin{enumerate}[(a)]
	\item $\Rank(H) = 2n + 1$ implies that $|H|$ is odd, and $|H|$ is $3/2$-transitive on $C$ with all $H_\alpha$-orbits in $C \setminus \{\alpha\}$ of size $(c-1)/(2n)$;
	\item $\Rank(K) = 2m + 1$ implies that $|K|$ is odd, and $K$ is $3/2$-transitive on $\Sigma$ with all $K_C$-orbits in $\C \setminus \{C\}$ of size $(d-1)/(2m)$.
	\end{enumerate}
\end{lemma}

\begin{proof}
Suppose first that $\Rank(H) = 2n + 1$. By Proposition \ref{prop:HK}, $(\Rank(H) - 1)/2 \leq |\OO_H| = \PairRank(H) \leq n$. Our assumption that $\Rank(H) = 2n + 1$ therefore implies that equality holds, and hence $\Sigma \neq \Sigma^*$ for each non-trivial $H$-orbital $\Sigma$. By a result in permutation groups (see, for instance, \cite[Lemma 2.27]{PS}), $H$ has odd order. Also $|\OO_H| = n$, and by Proposition~\ref{prop:HK}, each of the $n$ summands in \eqref{eq:n} is a positive integer. Hence each of these summands is equal to $1$, that is, $u_\Sigma = (c-1)/n$ for each non-trivial $H$-orbital $\Sigma$. Thus, for each such $\Sigma$, we have $|\Sigma(\alpha)| = u_\Sigma/2 = (c-1)/(2n)$ since $\Sigma \ne \Sigma^*$, and in particular $H$ is $3/2$-transitive on $C$. The proof of part (b) is similar.
\end{proof}

%%%%%%%%%%%%%%%%%%%%%%%%%%%%%%%%%%%%%%%%%%%%%%%%%%%%%%%%%%%%%%%%%%%

\section{Exploring examples: Proof of Theorem~\ref{theorem:exist}} \label{sec:existence}

In this section we explore the examples required to prove Theorem~\ref{theorem:exist}. For the first part we investigate the class of projective plane examples given in \cite[Example on p. 232]{OPP} and mentioned in \cite[p. 312]{P01}. We show that these examples satisfy all the conditions in Lemma \ref{lemma:maxrank} parts (a) and (b), and that in this family there are designs for which the Delandtsheer--Doyen parameters are (simultaneously) arbitrarily large. This therefore will prove Theorem~\ref{theorem:exist}(a).

\begin{example} \label{ex:projpl}
Let $q$ be a prime power such that $q^2 + q + 1$ is not prime, say $q^2 + q + 1 = cd$ where $c, d \geq 2$. Let $\D = (\PP,\B)$ where $\PP$ and $\B$ are the points and lines, respectively, of the Desarguesian projective plane $\PG_2(q)$. Then $v = |\PP| = q^2 + q + 1 = cd$ and $k = q+1$. Let $G$ be  a \emph{Singer cycle}, that is, a cyclic subgroup of automorphisms of $\D$ of order $cd$ acting regularly on $\PP$. Then $G$ is also transitive (in fact regular) on $\B$ (see \cite[Result 1 of Section 2.3]{Dem}). Also $G$ preserves a partition $\C$ of $\PP$ into $d$ classes of size $c$,
namely the set of orbits in $\PP$ of the unique cyclic (normal) subgroup of $G$ of order $c$. 
\end{example}

The information given in Example~\ref{ex:projpl} tells us that $\mathscr{D}$ is a $2$-$(cd, q+1, 1)$ design admitting $G = \mathbb{Z}_{cd}$ as a point-imprimitive, block-transitive group of automorphisms. We now prove the other assertions mentioned above.

\begin{lemma}\label{lemma:projpl}
Let $\D$, $G$, $\C$ be as in Example~$\ref{ex:projpl}$. 
	\begin{enumerate}[(a)]
	\item The Delandtsheer--Doyen parameters $(m,n)$ for $\D$ relative to $G$ and $\C$ are 
	\[ m = \frac{d-1}{2} \quad \mbox{and} \quad n = \frac{c-1}{2}. 	\]
	\item The group $G$ is permutationally isomorphic to a subgroup of $H \wr K$ where $H = \mathbb{Z}_c$, the group induced on a class of $\C$, and $K = \mathbb{Z}_d$, the group induced on $\C$.
	\item $\Rank(H) = 2n+1$ and $\Rank(K) = 2m+1$, the upper bounds of Theorem~$\ref{theorem:rankbounds}$.
	\item For any $N > 0$ there exists $q$ such that $q^2+q+1 = cd$ with both $m > N$ and $n > N$.
	\end{enumerate}
\end{lemma}

\begin{proof}
Since $G = \mathbb{Z}_{cd}$, the group induced on each class is $H = \mathbb{Z}_c$, and the group $G$ induced on $\C$ is $K = \mathbb{Z}_d$. By \cite[Theorem 5.5]{PS}, $G$ is permutationally isomorphic to a subgroup of $H \wr K$. Moreover, $H$ and $K$ are regular of degree $c$ and $d$, respectively, and hence in particular each is $3/2$-transitive. Thus $\Rank(H) = c$ and $\Rank(K) = d$, and by Theorem~\ref{theorem:DD}, the Delandtsheer--Doyen parameters $(m,n)$ are such that $c = \left(\binom{k}{2} - n\right)/m$ and $d = \left(\binom{k}{2} - m\right)/n$.  By Lemma \ref{counting}(a) $\binom{k}{2}= n \cdot \frac{cd-1}{c-1} = n \cdot \frac{q^2 + q}{c-1}.$
%From these expressions for $c$ and $d$, we get $\binom{k}{2} = mc + n = nd + m$, so that $m(c-1) = n(d-1)$. Hence
%	\[ \binom{k}{2} = mc + n = \frac{cn(d-1)}{c-1} + n = n \cdot \frac{c(d-1) + (c-1)}{c-1} = n \cdot \frac{cd-1}{c-1} = n \cdot \frac{q^2 + q}{c-1}. \]
However, also $\binom{k}{2} = q(q + 1)/2$ since $k=q+1$, and therefore $c = 2n + 1$. From Lemma \ref{counting}(b) we obtain $d = m(c-1)/n + 1 = 2m + 1$. Thus $\Rank(H) = 2n + 1$ and $\Rank(K) = 2m + 1$, which are the maximum possible values by Proposition~\ref{prop:HK}, and are the upper bounds of Theorem~\ref{theorem:rankbounds}. This proves parts (a)--(c).

To show that $n$ and $m$ can simultaneously be arbitrarily large, consider $q = p^{2f}$ for a prime $p$ and an integer $f$. Then
	\[ cd = q^2 + q + 1 = \frac{q^3 - 1}{q-1} = \frac{p^{3f} - 1}{p^f - 1} \cdot \frac{p^{3f} + 1}{p^f + 1} = \left(p^{2f} + p^f + 1\right)\left(p^{2f} - p^f + 1\right). \]
We may take $d = p^{2f} + p^f + 1$ and $c = p^{2f} - p^f + 1$, so that by part (a), $m = p^f\left(p^f + 1\right)/2$ and $n = p^f\left(p^f - 1\right)/2$. For any fixed prime $p$ and any given bound $N$, we can find $f$ such that $p^f\left(p^f - 1\right)/2 > N$, and then we will have $m, n > N$, so part (d) holds.
\end{proof}

Theorem~\ref{theorem:exist}(a) follows from Lemma~\ref{lemma:projpl}. Now we establish the other parts of Theorem~\ref{theorem:exist} using results from \cite{CP, NNOPP, OPP}. Namely, Lemma~\ref{lemma:existbc}(a) proves Theorem~\ref{theorem:exist}(b), and Lemma~\ref{lemma:existbc}(b) proves Theorem~\ref{theorem:exist}(c). Note that $m = n = 1$ implies that $c = d = \binom{k}{2} - 1$ by Theorem~\ref{theorem:DD}.

\begin{lemma} \label{lemma:existbc}
Let $k$ be an integer, $k \geq 3$, and let $c = \binom{k}{2} - 1$.
	\begin{enumerate}[(a)]
	\item There exist at least three pairwise non-isomorphic $2$-$(c^2, k, \la)$  designs (for some values of $\la$), each admitting $H \wr K$ as a block-transitive, point-imprimitive group of automorphisms with Delandtsheer--Doyen parameters $(m,n) = (1,1)$, where $H = K = \Sym(c)$ and hence with $\Rank(X) = \PairRank(X) + 1 = 2$, for $X \in \{H,K\}$.
	\item There exists a $2$-$(c^2, k, \la)$ design $\D$ (for some $\la$) as in Theorem~$\ref{theorem:rankbounds}$, with $n = m = 1$ and $\max\{\Rank(H), \Rank(K)\} = 3$, if and only if $k = 8$ and $c = 27$. Moreover, in the case $\la = 1$, there are up to isomorphism exactly 
	$467$ such designs.
	\end{enumerate}
\end{lemma}

\begin{proof}
(a) For each $k \geq 3$, constructions for three pairwise non-isomorphic designs with these parameters are described in the text after \cite[Theorem 5.1]{CP}, and for each design the subgroups $H = K = \Sym(c)$ satisfy $\Rank(H) = \Rank(K) = 2$ and $\PairRank(H) = \PairRank(K) = m = n = 1$.

(b) Suppose that, for some $\lambda$, there exists a $2$-$(c^2, k, \la)$ design  as in Theorem~\ref{theorem:rankbounds} with $n = m = 1$ and $\max\{\Rank(H), \Rank(K)\} = 3$. Then it follows from \cite[Theorem 5.2]{CP} that $k = 3, 4, 5$, or $8$ (because otherwise $\Rank(H) = \Rank(K) = 2$), and so $c = 2, 5, 9$, or $27$, respectively. By assumption some $X \in \{H,K\}$ has rank $3$, and hence $\PairRank(X) = 1$ by Theorem~\ref{theorem:rankbounds}, that is to say, $X$ is transitive on unordered pairs. However there is no transitive rank $3$ group of degree $c = 2$, and there is no transitive rank $3$ group of degree $c = 5$ or $c = 9$ that is transitive on unordered pairs (by \cite[Proposition 3.1]{Kan}, or see \cite[Theorem 9.4B]{DM}). Hence $k = 8$ and $c = 27$. In \cite{NNOPP,OPP}, it is proved that there are, up to isomorphism, exactly $467$ examples of $2$-$(729,8,1)$ designs (linear spaces) with these properties.
\end{proof}

%%%%%%%%%%%%%%%%%%%%%%%%%%%%%%%%%%%%%%%%%%%%%%%%%%%%%%%%%%%%%%%%%%%

\section{New design construction} \label{sec:construction}

In this section, we will construct block-transitive imprimitive designs with $m=1$ such that $\Rank(H) = \PairRank(H) + 1 = n + 1$ and $\Rank(K) = \PairRank(K) + 1 = m + 1 = 2$, for some fixed values of $n$ and of $c$.

Let $\F$ be a field of order $c=p^a$ such that $c \equiv 1 \pmod{2n}$, and let $\zeta$ be a primitive element of $\F$. Let $H = N \rtimes \lan \zeta^n \ra$ be the subgroup of the affine group $\AGL(1,c)$ acting on $\F$, where $N$ is the group of translations, and we identify $\zeta^n$ with multiplication by $\zeta^n$. Note that $\lan \zeta^n \ra$ contains $-1$ since $c \equiv 1 \pmod{2n}$. We record some information about the permutation action of $H$ on $\F$. The assertions are straightforward to check and details are left to the reader.

\begin{lemma} \label{lemma:H}
Let $\F$, $c$, $\zeta$, $N$, and $H$ be as above. Then the following hold:
	\begin{enumerate}[(a)]
	\item $N$ is regular on $\F$, and $H_0 = \lan \zeta^n \ra$ is the stabiliser in $H$ of $0 \in \F$.
	\item The $H_0$-orbits in $\F \setminus \{0\}$ are $\Delta_i(0) = \left\{ \zeta^{i+jn} \mid 0 \leq j < (c-1)/n \right\}$ with associated $H$-orbital $\Delta_i = (0,\zeta^i)^H$, for $0 \leq i < n$. Each $\Delta_i$ is self-paired (since $-1 \in H_0$).
	\item The $H$-orbits on $2$-subsets of $\F$ are $\OO_i = \{ \{\alpha,\beta\} \mid (\alpha,\beta) \in \Delta_i \}$, for $0 \leq i < n$,   each of which has size $c(c-1)/(2n)$. Moreover, the setwise stabiliser of each pair of points is cyclic of order $2$.
	\item $\Rank(H) = \PairRank(H) + 1 = n + 1$
	\end{enumerate}
\end{lemma}

We use this group $H$ in the design construction. The point-imprimitive group of automorphisms will be $G := H \wr K$, where $K = \Sym(d)$ is the symmetric group on $R = \mathbb{Z}_d$. Note that $G$ acts imprimitively on $\F \times R$ as follows (see \cite[Lemma 5.4]{PS}):
	\[ (x,j)^{(h_1,\dots,h_d)\sigma} = \left(x^{h_j},j^\sigma\right), \quad \mbox{for $(h_1,\dots,h_d) \in H^d$, $\sigma \in K$, and $(x,j) \in \F \times R$}. \]
The following properties of this $G$-action are again straightforward to check and details are left to the reader.

\begin{lemma} \label{lemma:G}
Let $\F$, $R$, and $G$ be as above. Then the following hold:
	\begin{enumerate}[(a)]
	\item The partition $\C$ of $\F \times R$, with classes $C_j = \{ (x,j) \mid x \in \F \}$ for $j \in R$, is nontrivial and $G$-invariant.
	\item $G$ is transitive on the set $\OO_{out}$ of $\C$-outer pairs from $\F \times R$, and $|\OO_{out}| = c^2d(d-1)/2$.
	\item $G$ has exactly $n$ orbits on the set of $\C$-inner pairs from $\F \times R$, namely 
	\[
	\OO_{inn,i} = \{ \{(x,j), (y,j)\} \mid \{x,y\} \in \OO_i, j \in R \},\quad
	 \mbox{ for $0 \leq i < n$}. 
	\]
Moreover, for each $i$ we have $|\OO_{inn,i}| = dc(c-1)/(2n)$.
	\end{enumerate}
\end{lemma}
 
For the construction below to work, we need some conditions on $n$ and $c$ which are exactly the conditions described in  Definition~\ref{def:nc}. 
Suppose now that $[n,c]$ is a useful pair, as in Definition~\ref{def:nc}. Then $c = p^a$, for some odd prime $p$ and $a \geq 1$, and $c \equiv 1 \pmod{2n}$, $c + n = \binom{k}{2}$ for some integer $k \geq 2n$, and $n\geq 2$. First we derive an upper bound and an improved lower bound for $k$.

\begin{lemma} \label{lemma:kbounds}
Let $n$, $c$, and $k$ be as above and let $d = 1 + (c-1)/n$. Then $2n + 2 \leq k \leq n + d$.
\end{lemma}

\begin{proof}
From Definition~\ref{def:nc} we have $c + n = k(k-1)/2 \geq n(k-1)$.
Suppose first, for a contradiction, that $k > n + d$. Then
	\[ 
	k > n + d = n + 1 + \frac{c-1}{n} = n - \frac{1}{n} + \frac{c+n}{n} \geq  n - \frac{1}{n} + (k-1), 
	\]
and hence $n < 1 + 1/n$ which is not possible for any $n \geq 2$. Hence $k \leq n + d$. Finally, if $k = 2n$ or $k = 2n + 1$, then
	\[ c = \binom{k}{2} - n \in \{ 2n(n-1), 2n^2 \}, \]
and in either case $c \equiv 0 \pmod{2n}$, which is a contradiction. Hence $k \geq 2n + 2$.
\end{proof}

\begin{construction} \label{con:newdesign}
Let $[n,c]$ be a useful pair as in Definition~\ref{def:nc}, with $c = p^a$ for some odd prime $p$ and $a \geq 1$, $d = 1 + (c-1)/n$, and $k$, $\F$, $R$, $G$, and $\C$ as above. Define the design $\D = (\PP,\B)$ to have point set $\PP := \F \times R$ and block set $\B := B^G$, where $B \subseteq \PP$ is given by
	\begin{equation}\label{eqb}
	B = \{ (0,i), (\zeta^i,i) \mid 0 \leq i \leq n-1 \} \cup \{ (0,i) \mid n \leq i \leq k-n-1 \}. 	
\end{equation}	 
Note that $n+1 \leq k-n-1 \leq d-1$, by Lemma~\ref{lemma:kbounds}, so the second set in the union defining $B$ has size  $k-2n\geq 2$, and $B$ is a well defined $k$-subset of $\PP$.
\end{construction}

\begin{theorem} \label{theorem:newdesign}
Let $\D$ and $G$ be as in Construction~$\ref{con:newdesign}$. Then $\D = (\PP,\B)$ is a $2$-$(cd, k, \la)$ design (for some $\la$) and $G = H \wr K$ is a block-transitive, point-imprimitive group of automorphisms leaving invariant the point-partition $\C$. Moreover:
	\begin{enumerate}[(a)]
	\item the Delandtsheer--Doyen parameters of $\D$ are $(m,n) = (1,n)$; and
	\item $\Rank(H) = \PairRank(H) + 1 = n + 1$ and $\Rank(K) = \PairRank(K) + 1 = m + 1 = 2$.
	\end{enumerate}
\end{theorem}

\begin{proof}
By \eqref{eqb} and Lemma~\ref{lemma:G}, for each $i=0,\dots,n-1$, the set $B$ contains exactly $n_i = 1$ inner pair from $\OO_{inn,i}$, namely $\{ (0,i), (\zeta^i,i) \}$. Hence $B$ contains exactly $n_{out} := \binom{k}{2} - n$ outer pairs. By \cite[Proposition 1.3]{CP}, $\D$ is a $2$-design if and only if, for any $G$-orbit $\OO$ in the set of $2$-subsets of $\PP$, the ratio $|\{ \{x,y\} \mid \{x,y\} \in \OO, \ \{x,y\} \subseteq B \}|/|\OO|$ is independent of $\OO$, that is to say, $\D$ is a $2$-design if and only if
	\[ 
	\frac{n_0}{|\OO_{inn,0}|} = \dots = \frac{n_{n-1}}{|\OO_{inn,n-1}|} = \frac{n_{out}}{|\OO_{out}|}. 
	\]
Hence, by Lemma~\ref{lemma:G}, $\D$ is a $2$-design if and only if
	\[ 
	\frac{1}{dc(c-1)/(2n)} = \frac{\binom{k}{2}-n}{c^2d(d-1)/2}, 
	\]
or equivalently,
	\[ 
	\binom{k}{2}-n = \frac{nc(d-1)}{c-1}. 
	\]
From the definition of $d$, this is equivalent to $\binom{k}{2} - n = c$, and this equality holds by Definition~\ref{def:nc}. Hence $\D$ is a $2$-design with $cd$ points and block size $k$, that is to say, a $2$-$(cd, k, \la)$ design for some $\la$. By Definition~\ref{def:nc} and Lemma~\ref{lemma:kbounds}, the class size $c = \binom{k}{2} - n$ and the number of classes $d = (n + c - 1)/n = \left(\binom{k}{2} - 1\right)/n$, and it follows from Theorem~\ref{theorem:DD} that the Delandtsheer--Doyen parameters are $(m,n) = (1,n)$. This proves part (a). The assertions in part (b) follow from Lemma~\ref{lemma:H}(d) for $H$,  and the fact that $K = \Sym(d)$ has $\Rank(K) = \PairRank(K) + 1 = 2 = m+1$.
\end{proof}

\begin{remark} \label{rem:newdesign-la}
It is possible, but not very insightful, to determine the value of $\la$ in Theorem~\ref{theorem:newdesign}. However, for completeness we do so here. 
%We count the number of pairs $(\rho,B')$ where $\rho \in \OO_{inn,0}$, $B' \in \B$, and $\rho \subseteq B'$. This yields $|\OO_{inn,0}|\la = |\B|n_0$, and evaluating this we find, using Lemma~\ref{lemma:G}, $\big(dc(c-1)/(2n)\big) \cdot \la = |\B|\cdot 1$, that is,

 By Lemma \ref{counting}(c) and since $m=1$, we have 
	\( \la =% \frac{2n|\B|}{dc(c-1)} =
	\frac{2|\B|}{cd(d-1)}. \)
Recall (from the definition of $H$) that a point stabiliser in $H$ is cyclic of order $(c-1)/n = d-1$, and (from Lemma~\ref{lemma:H}(c)) that the stabiliser in $H$ of each unordered pair of points is cyclic of order $2$. From these and the definition of $B$ we see that
	\[ G_B = (\mathbb{Z}_2\wr \Sym(n)) \times (\mathbb{Z}_{d-1} \wr \Sym(k-2n)) \times (H \wr \Sym(d-k+n)). \]
Hence
	\begin{align*}
	|\B| = |G:G_B|
	&= \frac{|H|^d d!}{2^n(d-1)^{k-2n} |H|^{d-k+n} n! (k-2n)! (d-k+n)!} \\
	&= \frac{c^{k-n} (d-1)^n}{2^n} \cdot \frac{d!}{n!(k-2n)!(d-k+n)!}
	\end{align*}
which yields
	\begin{align*}
	\lambda
	&= \frac{2}{cd(d-1)} \cdot c^{k-n} \left(\frac{d-1}{2}\right)^n \frac{d!}{n! (k-2n)! (d-k+n)!} \\
	&= c^{k-n-1} \left(\frac{d-1}{2}\right)^{n-1} \frac{(d-1)!}{n!(k-2n)!(d-k+n)!}.
	\end{align*}
For example, the smallest useful pair is $[n,c] = [2,13]$, and for this pair the value of $\la$ is $197730$.
\end{remark}

Note there are many useful pairs, see Table \ref{table:nc}, but they do not exist for every $n$, as proved in the following lemma. We comment on the existence of useful pairs in Subsection \ref{sub:useful}.

\begin{lemma} \label{lemma:nwithnousefulpair}
If $[n,c]$ is a useful pair, then $n\notin\{6,10,15\}$. Moreover, if $n\in\{6,10,15\}$ and $[n,c]$ satisfies all the conditions of a useful pair except  that $k< 2n$, then  $[n,c,k,d]$ is one of $[6, 49, 11, 9]$, $[10, 81, 14, 9]$, or $[15, 121, 17, 9]$.
\end{lemma}
\begin{proof} Let $n\in\{6,10,15\}$, and let  $[n,c]$ satisfy all the conditions for a useful pair except possibly $k\geq 2n$. Then $ \binom{k}{2}-n\equiv 1\pmod {2n}$, so $k(k-1)\equiv 2n+2\pmod {4n}$. 
For $n=6$ this implies $k\equiv 11$ or $14  \pmod{4n}$, for $n=10$ this implies $k\equiv 14,19,22$ or $27  \pmod{4n}$, and for  $n=15$ this implies $k\equiv 17,29,32$ or $44  \pmod{4n}$.
So $k=4nb+r$ for some non-negative integer $b$ and some residue $r$ as in the previous sentence, and $\binom{k}{2}-n$, as a quadratic polynomial $g(b)$, happens to factorise as a product of two linear polynomials with integer coefficients as in Table~\ref{tab:withnousefulpair}.

\begin{table}[h]
    \centering
\begin{tabular}{l|l|l}
\hline
$n$&$k$&$g(b)= \binom{k}{2}-n$\\
\hline
$6$&$24b+11$&$(12b+7)(24b+7)$\\
&$24b+14$&$(12b+5)(24b+17)$\\
\hline
$10$&$40b+14$&$(20b+9)(40b+9)$\\
&$40b+19$&$(20b+7)(40b+23)$\\
&$40b+22$&$(20b+13)(40b+17)$\\
&$40b+27$&$(20b+11)(40b+31)$\\
\hline
$15$&$60b+17$&$(30b+11)(60b+11)$\\
&$60b+29$&$(30b+17)(60b+23)$\\
&$60b+32$&$(30b+13)(60b+37)$\\
&$60b+44$&$(30b+19)(60b+49)$\\
\hline
\end{tabular}
    \caption{Factorisations of $\binom{k}{2}-n$ for the proof of Lemma~\ref{lemma:nwithnousefulpair}}
    \label{tab:withnousefulpair}
\end{table}

In each case the polynomial $g(b)$ has the form $(2nb+x)(4nb+y)$ for some integers $x, y$ satisfying $1<x\leq y$.  Since we must have $g(b)=p^a$ for some prime $p$ and positive integer $a$, it follows that $2nb+x=p^e$ and $4nb+y=p^f$ for integers $e, f$ such that $0<e\leq f$ and $e+f=a$. Therefore $y-2x=p^f-2p^e=p^e(p^{f-e}-2)$. 

For $n=6$, $y-2x=p^e(p^{f-e}-2)=\pm 7$, so $p=7$, $e=f=1$, and hence $x=y=7$ which forces  $b=0$. Thus  $[n,c,k,d]=[6, 49, 11, 9]$. These values do not satisfy the condition $k\geq 2n$ so $[6,49]$ is not a useful pair.

For $n=10$, $y-2x=p^e(p^{f-e}-2)=\pm 9$, so $p=3$, $e=2$ and $f=2$ or $3$. If $f=2$, then $b=0$ and $x=y=9$, so  $[n,c,k,d]=[10, 81, 14, 9]$. These values do  not satisfy the condition $k\geq 2n$ so $[10,81]$ is not a useful pair. If $f=3$, then $20b+x=9,40b+y=27$, which force $b=0$ and $y=27$. However there is no line in Table~\ref{tab:withnousefulpair} with $y=27$.

For $n=15$, $y-2x=p^e(p^{f-e}-2)=\pm 11$, so $p=11$, $e=f=1$. This forces $b=0$ and $x=y=11$, so  $[n,c,k,d]=[15, 121, 17, 9]$. Again these values do not satisfy the condition $k\geq 2n$ so $[15,121]$ is not a useful pair.
\end{proof}

\end{document}